\documentclass[12pt,draft]{amsart}
\usepackage{amsmath,amssymb,amsthm,bm}
\usepackage{graphicx,enumerate}
\usepackage{layout}

\theoremstyle{plain}
\newtheorem{theorem}{Theorem}[section]

\theoremstyle{definition}

\numberwithin{equation}{section}

\title[A rationality problem of some Cremona transformation]
{\Large A rationality problem of \vspace*{2mm}\\some Cremona transformation}

\subjclass[2000]{Primary 14E07, 14E08, 13A50, 12F20.}

\keywords{Rationality problem, Cremona transformations, linear actions, monomial group actions.}

\begin{document}
\maketitle
\begin{center}
\begin{tabular}{lll}
Akinari Hoshi & &Ming-chang Kang\\
Department of Mathematics & &Department of Mathematics\\
School of Education & and\hspace*{3mm} &National Taiwan University\\
Waseda University& &Taipei, Taiwan \\
Tokyo, Japan & &E-mail: \texttt{kang@math.ntu.edu.tw}\\
E-mail: \texttt{hoshi@ruri.waseda.jp} &
\end{tabular}
\end{center}
\vspace*{5mm}

\noindent Abstract.
Let $k$ be any field, $k(x,y)$ be the rational function field of two variables over $k$.
Let $\sigma$ be a $k$-automorphism of $k(x,y)$ defined by
\begin{align*}
\displaystyle{\sigma(x)\,=\,\frac{-x(3x-9y-y^2)^3}{(27x+2x^2+9xy+2xy^2-y^3)^2},\quad
\sigma(y)\,=\,\frac{-(3x+y^2)(3x-9y-y^2)}{27x+2x^2+9xy+2xy^2-y^3}}.
\end{align*}
\noindent Theorem. The fixed field $k(x,y)^{\langle\sigma\rangle}$
is rational (= purely transcendental) over $k$. Embodied in the
proof of the above theorem are several general guidelines for
solving the rationality problem of Cremona transformations, which
may be applied elsewhere. \vspace*{2mm}

\section{Introduction}\label{sec-intro}

Let $k$ be any field, $k(x_1,\ldots,x_n)$ be the rational function field of $n$ variables.
(It is not necessary to assume that $k$ is algebraically closed.) By a Cremona transformation
on $\mathbb{P}^n$ we mean a $k$-automorphism $\sigma$ on $k(x_1,\ldots,x_n)$, i.e.
\begin{align}
\sigma\ :\ k(x_1,\ldots,x_n)\ \longrightarrow\ k(x_1,\ldots,x_n)\label{auto}
\end{align}
where $\sigma(x_i)\in k(x_1,\ldots,x_n)$ for each $1\leq i\leq n$ and $\sigma$
is an automorphism.
We will denote by $\mathrm{Cr}_n$ the group of all Cremona transformations on $\mathbb{P}^n$.
The purpose of this note is to consider whether $k(x_1,x_2)^G$ is rational (= purely
transcendental) over $k$ where $G$ is some finite subgroup of $\mathrm{Cr}_2$.

Note that, if $k$ is algebraically closed, then $k(x_1,x_2)^G$ is
rational over $k$ by Zariski-Castelnuovo's Theorem \cite{Za}. On
the other hand, if the group $G$ consists of automorphisms
$\sigma$ such that, in (\ref{auto}), $\sigma(x_i)$ are homogeneous
linear polynomials (resp. monomials) in $x_1,\ldots,x_n$, then the
group action of $G$ on $k(x_1,\ldots,x_n)$ is the usual linear
action (resp. the monomial group action). The rationality problem
of linear actions or the monomial group actions has been
investigated extensively. See, for examples, \cite{Sw, KP, HK1,
HK2, HR}. It seems that not many research works are devoted to the
rationality problem of ``genuine'' Cremona transformations, i.e.
the $\sigma(x_i)$ in (\ref{auto}) are, instead of linear
polynomials or monomials, but rational functions with total
degrees high enough, say, $\geq 4$. As far as we know, only ad hoc
techniques can be found in the literature for solving the
rationality problems of Cremona transformations.

The main result of this note is the following theorem.

\begin{theorem}\label{th1}
Let $k$ be any field and $k(x_1,x_2)$ be the rational function
field of two variables over $k$. Let $\sigma\in \mathrm{Cr}_2$
defined by

\begin{align*}
\sigma\ :\ k(x_1,x_2)\longrightarrow k(x_1,x_2)
\end{align*}
 where
\begin{align*}
\sigma(x_1)\,&=\,\frac{-x_1(3x_1-9x_2-x_2^2)^3}{(27x_1+2x_1^2+9x_1x_2+2x_1x_2^2-x_2^3)^2},\\
\sigma(x_2)\,&=\,\frac{-(3x_1+x_2^2)(3x_1-9x_2-x_2^2)}{27x_1+2x_1^2+9x_1x_2+2x_1x_2^2-x_2^3}.
\end{align*}
Then $k(x_1,x_2)^{\langle \sigma\rangle}\,:=\,\{ f\in k(x_1,x_2)\ : \sigma(f)\,=\,f\}$
is rational over $k$.
\end{theorem}

Note that $\sigma^2\,=\,1$.

The above theorem was given in \cite[Theorem 10]{HM}.
Unfortunately the proof in \cite{HM} contains a few mistakes. For
examples, the $\sigma_1, \sigma_2$ defined in \cite[p.25]{HM} are
not automorphisms. We will give another proof of Theorem \ref{th1}
in Section 2 (when $\mathrm{char}\ k\neq 2,3$) and Section 3 (when
$\mathrm{char}\ k= 2$, or $3$). Our proof is completely different
from that in \cite{HM}. We hope that this proof will be helpful to
people working on the rationality problem of Cremona
transformations, because it contains systematic methods for
attacking the rationality problem. (See Step 1, Step 2 and Step 5
of Section 2, in particular.) In keeping with the spirit of the
proof in Section 2 we give another proof of the case
$\mathrm{char}\, k\,=\,2$ and the case $\mathrm{char}\, k\,=\,3$
in Section 4 and Section 5 respectively.

Many rationality problems arise from the study of moduli spaces of
some geometric configurations. The rationality problem in Theorem
\ref{th1} arose in the study of the moduli of cubic generic
polynomials. See \cite{HM}.

Some symbolic computations in this note are carried out with the
aid of ``Mathematica'' \cite{Wo}.

Finally we will emphasize that it is unnecessary to assume that
the base field $k$ is algebraically closed or any restriction on
the characteristic of $k$.

\bigskip

\section{The case $\mathrm{char}\, k \neq 2,3$}\label{sec-neq23}

Throughout this section, we assume that $\mathrm{char}\ k\neq 2,3$. \\
%

Step $1$. Note that $\sigma$ induces a birational map on
$\mathbb{P}^2$. We will find some irreducible exceptional divisors
of this rational map. Clearly the curve defined by
$3x_1-9x_2-x_2^2\,=\,0\,$ is one of the candidates. Taking its
image $\sigma(3x_1-9x_2-x_2^2)$, we will find another polynomial.
Thus, define
\begin{align}
y_1\,&=\, 3x_1-9x_2-x_2^2,\nonumber\\
y_2\,&=\, 27x_1+9x_1x_2+x_2^3,\label{y123}\\
y_3\,&=\, -27x_1-2x_1^2-9x_1x_2-2x_1x_2^2+x_2^3.\nonumber
\end{align}
With the aid of computers, it is easy to see that
\begin{align}
\sigma\ :\ y_1\,\longmapsto\, y_1y_2^2y_3^{-2},\quad y_2\,\longmapsto\, y_1^3y_2^2y_3^{-3},\quad
y_3\,\longmapsto\, y_1^3y_2^3y_3^{-4}. \label{actsy}
\end{align}
Note that the determinant of the exponents of the above map is
\[
\mathrm{det}\,
\left(
\begin{array}{ccc} 1 & 3 & 3\\ 2 & 2 & 3 \\ -2~ & -3~ & -4~\end{array}
\right) \,=\, 1.
\]
Thus the action of $\sigma$ on $k(y_1,y_2,y_3)$ can be lifted to
$k(Y_1,Y_2,Y_3)$ ($Y_1,Y_2,Y_3$ are algebraically independent over
$k$) and induces a monomial action on $k(Y_1,Y_2,Y_3)$. But we
will not use this fact in the following steps.\\

Step $2$. Luckily we find that $k(y_1,y_2,y_3)\,=\,k(x_1,x_2)$. In
fact, from (\ref{y123}), we may eliminate $x_2$ and get two
polynomial equations of $x_1$ with coefficients in
$k(y_1,y_2,y_3)$; applying the Euclidean algorithm to these two
polynomials, we may show that $x_1\in k(y_1,y_2,y_3)$.

More explicitly, with the aid of computers, we will find (i) the
expressions of $x_1,x_2$ in terms of $y_1,y_2,y_3$, and (ii) a
polynomial equations of $y_1,y_2,y_3$. We get
\begin{align*}
x_1\,&=\, \frac{-2y_1^3-729y_2+27y_1y_2-2y_2^2-729y_3+27y_1y_3}{108(y_2+y_3)},\\
x_2\,&=\, \frac{-2y_1^4+9y_1^2y_2-2y_1y_2^2+9y_1^2y_3+81y_2y_3+81y_3^2}{18(y_1^3-y_2y_3)},
\end{align*}
\begin{align}
f(y_1,y_2,y_3)\,&=\,
2y_1^6+729y_1^3y_2-27y_1^4y_2+4y_1^3y_2^2-27y_1y_2^3+2y_2^4\label{fy1y2y3}\\
&+729y_1^3y_3-27y_1^4y_3-27y_1y_2^2y_3+729y_2y_3^2+729y_3^3\,=\,0.\nonumber
\end{align}

Step $3$.
The map of $\sigma$ defined in (\ref{actsy}) can be simplified as follows.
Define
\[
z_1\,=\,y_2^{-1}y_3,\quad z_2\,=\,y_1y_2^{-1},\quad z_3\,=\,y_1^{-2} y_3.
\]
It follows that $k(y_1,y_2,y_3)\,=\,k(z_1,z_2,z_3)$ and
\begin{align}
\sigma\ :\ z_1\,\longmapsto\, 1/z_1,\quad z_2\,\longmapsto\, z_3\,\longmapsto\, z_2.\label{actsz}
\end{align}
The relation $f(y_1,y_2,y_3)\,=\,0\,$ in (\ref{fy1y2y3}) becomes
\begin{align}
g(z_1,z_2,z_3)\,=\,2z_1^2z_2^2+4z_1z_2z_3-27z_1z_2^2z_3-27z_1^2z_2^2z_3+2z_3^2-27z_2z_3^2&
\label{gz1z2z3}\\
-27z_1z_2z_3^2+729z_2^3z_3^2+729z_1z_2^3z_3^2+729z_1z_2^2z_3^3+729z_1^2z_2^2z_3^3&\,=\,0.
\nonumber
\end{align}

Step $4$.
The map of $\sigma$ defined in (\ref{actsz}) is equivalent to
\begin{align*}
\sigma\ :\ z_2-z_3\,\longmapsto\, -(z_2-z_3),\quad \frac{1-z_1}{1+z_1}\,\longmapsto\,
-\frac{1-z_1}{1+z_1},\quad z_2+z_3\,\longmapsto\, z_2+z_3.
\end{align*}
Thus $k(x_1,x_2)^{\langle\sigma\rangle}\,=\,k(z_1,z_2,z_3)^{\langle\sigma\rangle}\,=\,k(u_1,u_2,u_3)$
where $u_1,u_2,u_3$ are defined by
\begin{align*}
u_1\,=\,(z_2-z_3)^2,\quad u_2\,=\,\Bigl(\frac{1-z_1}{1+z_1}\Bigr)\cdot (z_2-z_3),\quad u_3\,=\,z_2+z_3.
\end{align*}
The relation $g(z_1,z_2,z_3)\,=\,0\,$ in (\ref{gz1z2z3}) becomes
\begin{align}
108u_1u_2-729u_1^2u_2-16u_2^2-108u_1u_3-729u_1^2u_3+32u_2u_3-16u_3^2&\label{hu1u2u3}\\
-108u_2u_3^2+1458u_1u_2u_3^2+108u_3^3+1458u_1u_3^3-729u_2u_3^4-729u_3^5&\,=\,0.\nonumber
\end{align}

In conclusion, $k(x_1,x_2)^{\langle\sigma\rangle}$ is a field generated by $u_1,u_2,u_3$
over $k$ with the relation (\ref{hu1u2u3}).
We will simplify the relation (\ref{hu1u2u3}) to get two generators.\\

Step $5$.
The relation (\ref{hu1u2u3}) defines an algebraic surface.
However this algebraic surfaces contains singularities.
We will make some change of variables to simplify the singularities and the equation
(\ref{hu1u2u3}).
Define
\begin{align*}
v_1\,=\,u_1u_3^{-1},\quad v_2\,=\,u_2u_3^{-1},\quad v_3\,=\,u_3.
\end{align*}
Then $k(u_1,u_2,u_3)\,=\,k(v_1,v_2,v_3)$ and the relation (\ref{hu1u2u3}) becomes
\begin{align}
h(v_1,&v_2,v_3)\,=\, 16+108v_1-32v_2-108v_1v_2+16v_2^2\nonumber\\
&-108v_3+729v_1^2v_3+108v_2v_3+729v_1^2v_2v_3\label{hv1v2v3}\\
&-1458v_1v_3^2-1458v_1v_2v_3^2+729v_3^3+729v_2v_3^3\,=\,0.\nonumber
\end{align}
We will determine the singularities of $h(v_1,v_2,v_3)=0$ by
solving
\[
h(v_1,v_2,v_3)\,=\,\frac{\partial h}{\partial
v_1}(v_1,v_2,v_3)\,=\, \frac{\partial h}{\partial
v_2}(v_1,v_2,v_3)\,=\,0.
\]
We get $v_2-1\,=\,v_1-v_3\,=\,0$. Define
\[
w_1\,=\,v_1-v_3,\quad w_2\,=\,v_2-1,\quad w_3\,=\,v_3.
\]
Then $k(v_1,v_2,v_3)\,=\,k(w_1,w_2,w_3)$ and the relation (\ref{hv1v2v3}) becomes
\[
108w_1w_2 - 16w_2^2 - 1458w_1^2w_3 - 729w_1^2w_2w_3\,=\,0.
\]
The above equation is a linear equation in $w_3$.
Thus $w_3\in k(w_1,w_2)$.
It follows $k(w_1,w_2,w_3)\,=\,k(w_1,w_2)$.
We conclude that $k(x_1,x_2)^{\langle\sigma\rangle}\,=\,k(w_1,w_2,w_3)\,=\,k(w_1,w_2)$
is rational over $k$. \\

Step $6$.
We will give explicit formulae of $w_1,w_2$ in terms of $x_1,x_2$.
It is not difficult to find that
\begin{align*}
w_1\,&=\,\frac{-4(3x_1-9x_2-x_2^2)(27x_1+2x_1^2+9x_1x_2+2x_1x_2^2-x_2^3)}
{(27+x_1+9x_2+x_2^2)(27x_1^2+18x_1^2x_2-27x_1x_2^2+27x_2^3+2x_1x_2^3)},\\
w_2\,&=\,\frac{27(27x_1+2x_1^2+9x_1x_2+2x_1x_2^2-x_2^3)}
{27x_1^2+18x_1^2x_2-27x_1x_2^2+27x_2^3+2x_1x_2^3}.
\end{align*}
We also see
\begin{align*}
\frac{w_1}{w_2}\,=\,\frac{-4(3x_1 - 9x_2 - x_2^2)}{27(27 + x_1 + 9x_2 + x_2^2)}.
\end{align*}
Finally we obtain
\begin{align*}
k(x_1,x_2)^{\langle\sigma\rangle}\,=\,k\Bigl(\frac{3x_1 - 9x_2 - x_2^2}{27 + x_1 + 9x_2 + x_2^2},
\frac{27x_1+2x_1^2+9x_1x_2+2x_1x_2^2-x_2^3}{27x_1^2+18x_1^2x_2-27x_1x_2^2+27x_2^3+2x_1x_2^3}\Bigr).
\end{align*}

\bigskip

\section{The remaining cases}\label{sec-remaining}

Step 1. In this step, we assume that $\mathrm{char}\, k\,=\,2$.
Note that the automorphism $\sigma$ becomes
\begin{align*}
x_1\,\longmapsto\,
\frac{x_1(x_1+x_2+x_2^2)^3}{(x_1^2+x_1x_2+x_2^3)^2},\quad
x_2\,\longmapsto\,
\frac{(x_1+x_2^2)(x_1+x_2+x_2^2)}{x_1^2+x_1x_2+x_2^3}.
\end{align*}
Define
\begin{align*}
y_1\,=\,x_1+x_2+x_2^2,\quad y_2\,=\,x_2.
\end{align*}
Then we have $k(x_1,x_2)\,=\,k(y_1,y_2)$ and
\[
\sigma\ :\ y_1\,\longmapsto\, y_1,\quad y_2\,\longmapsto\,
\frac{y_1(y_1+y_2)}{y_1+y_2+y_1y_2}.
\]
Also define
\[
z_1\,=\,y_1,\quad z_2\,=\,\frac{y_1+y_2}{y_2}.
\]
It follows that $k(y_1,y_2)\,=\,k(z_1,z_2)$ and
\[
\sigma\ :\ z_1\,\longmapsto\, z_1,\quad z_2\,\longmapsto\, z_1
z_2^{-1}.
\]
Therefore we obtain
\[
k(x_1,x_2)^{\langle\sigma\rangle}\,=\,k(z_1,z_2)^{\langle\sigma\rangle}
\,=\,k\Bigl(z_1,z_2+\frac{z_1}{z_2}\Bigr)\,=\,k\Bigl(x_1+x_2+x_2^2,
\frac{x_1^2+x_1x_2^2+x_2^3}{x_2(x_1+x_2^2)}\Bigr).
\]

\bigskip
Step 2. In this step, we assume that $\mathrm{char}\, k\,=\,3$.
Note that the automorphism $\sigma$ becomes
\begin{align*}
x_1\,\longmapsto\, \frac{x_1x_2^6}{(x_1^2+x_1x_2^2+x_2^3)^2},\quad
x_2\,\longmapsto\, \frac{-x_2^4}{x_1^2+x_1x_2^2+x_2^3}.
\end{align*}
Define
\begin{align*}
y_1\,=\,x_1 x_2^{-2},\quad y_2\,=\,x_2^{-1}.
\end{align*}
It follows that $k(x_1,x_2)\,=\,k(y_1,y_2)$ and
\[
\sigma\ :\ y_1\,\longmapsto\, y_1,\quad y_2\,\longmapsto\,
-y_2-y_1-y_1^2.
\]
Hence we get
\[
k(x_1,x_2)^{\langle\sigma\rangle}\,=\,k(y_1,y_2)^{\langle\sigma\rangle}
\,=\,k\bigl(y_1,y_2(y_2+y_1+y_1^2)\bigr)\,=\,k\Bigl(\frac{x_1}{x_2^2},
\frac{x_1^2+x_1x_2^2+x_2^3}{x_2^5}\Bigr).
\]

\bigskip

\section{The case $\mathrm{char}\, k=2$}\label{sec-char2}

In this section, we assume that $\mathrm{char}\, k\,=\,2$. Recall
that the automorphism $\sigma$ is
\begin{align*}
x_1\,\longmapsto\, \frac{x_1(x_1+x_2+x_2^2)^3}{(x_1^2+x_1x_2+x_2^3)^2},\quad
x_2\,\longmapsto\, \frac{(x_1+x_2^2)(x_1+x_2+x_2^2)}{x_1^2+x_1x_2+x_2^3}.
\end{align*}
Define
\begin{align}
y_1\,=\, x_1,\quad y_2\,=\, x_1+x_2+x_2^2,\quad y_3\,=\, x_1+x_1x_2+x_2^3.
\end{align}
With the aid of computers, it is easy to see that
\begin{align*}
\sigma\ :\ y_1\,\longmapsto\, y_1y_2^3y_3^{-2},\quad y_2\,\longmapsto\, y_2,\quad
y_3\,\longmapsto\, y_2^3y_3^{-1}.
\end{align*}
From (4.1) , we find that
\begin{align}
x_2\,=\,\frac{y_2+y_3}{1+y_2} .
\end{align}
And therefore we have that $k(y_1,y_2,y_3)\,=\,k(x_1,x_2)$. Using
(4.1) to eliminate $x_1,x_2$, we obtain the relation
\begin{align}
f(y_1,y_2,y_3)\,=\,y_1+y_1y_2^2+y_2^3+y_3+y_2y_3+y_3^2\,=\,0.\label{relych2}
\end{align}
Define
\[
z_1\,=\,y_1^{-1}y_3,\quad z_2\,=\,y_2^2y_3^{-1},\quad z_3\,=\,y_2^{-1} y_3.
\]
It follows that $k(y_1,y_2,y_3)\,=\,k(z_1,z_2,z_3)$ and
\begin{align*}
\sigma\ :\ z_1\,\longmapsto\, z_1,\quad z_2\,\longmapsto\, z_3\,\longmapsto\, z_2.
\end{align*}
We find that the relation $f(y_1,y_2,y_3)\,=\,0\,$ in
(\ref{relych2}) becomes
\begin{align}
g(z_1,z_2,z_3)\,=\,1+z_1+z_2z_3+z_2^2z_3+z_2z_3^2+z_1z_2^2z_3^2\,=\,0.\label{relzch2}
\end{align}
Define
\begin{align*}
u_1\,=\,z_1,\quad u_2\,=\,z_2z_3,\quad u_3\,=\,z_2+z_3.
\end{align*}
Then we have
$k(x_1,x_2)^{\langle\sigma\rangle}\,=\,k(z_1,z_2,z_3)^{\langle\sigma\rangle}\,=\,
k(u_1,u_2,u_3)$ and the relation in (\ref{relzch2}) becomes
\begin{align}
1 + u_1 + u_2 + u_1u_2^2 + u_2u_3\,=\,0.\nonumber
\end{align}
Thus $u_3\in k(u_1,u_2)$. It follows that
$k(x_1,x_2)^{\langle\sigma\rangle}\,=\,k(u_1,u_2,u_3)\,
=\,k(u_1,u_2)$ is rational over $k$. It is easy to obtain the
formulae of the generators $u_1,u_2$ of
$k(x_1,x_2)^{\langle\sigma\rangle}$ in terms of $x_1,x_2$. Indeed
we have
\[
u_1\,=\, \frac{x_1}{x_1+x_1 x_2+x_2^3},\quad u_2\,=\,x_1+x_2+x_2^2.
\]

\bigskip

\section{The case $\mathrm{char}\, k=3$}\label{sec-char3}

In this section, we assume that $\mathrm{char}\, k\,=\,3$. Recall
that the automorphism $\sigma$ is
\begin{align*}
x_1\,\longmapsto\, \frac{x_1x_2^6}{(x_1^2+x_1x_2^2+x_2^3)^2},\quad
x_2\,\longmapsto\, \frac{-x_2^4}{x_1^2+x_1x_2^2+x_2^3}.
\end{align*}
Define
\begin{align*}
y_1\,=\,x_1,\quad y_2\,=\,-x_2,\quad y_3\,=\, x_1^2+x_1x_2^2+x_2^3.
\end{align*}
It is clear that $k(x_1,x_2)\,=\,k(y_1,y_2,y_3)$ and
\[
\sigma\ :\ y_1\,\longmapsto\, y_1y_2^6y_3^{-2},\quad y_2\,\longmapsto\, y_2^4y_3^{-1},\quad
y_3\,\longmapsto\, y_2^{15}y_3^{-4}.
\]
The map of $\sigma$ above can be simplified as follows.
Define
\[
z_1\,=\,y_1y_2^{-2},\quad z_2\,=\,y_2^{-4}y_3,\quad z_3\,=\,y_2^{-1}.
\]
It follows that $k(y_1,y_2,y_3)\,=\,k(z_1,z_2,z_3)$ and
\begin{align*}
\sigma\ :\ z_1\,\longmapsto\, z_1,\quad z_2\,\longmapsto\, z_3\,\longmapsto\, z_2.
\end{align*}
We also obtain the relation
\begin{align}
g(z_1,z_2,z_3)\,=\,z_1+z_1^2-z_2-z_3\,=\,0.\label{relzch3}
\end{align}
Thus $k(x_1,x_2)^{\langle\sigma\rangle}\,=\,k(z_1,z_2,z_3)^{\langle\sigma\rangle}\,=\,
k(u_1,u_2,u_3)$ where $u_1,u_2,u_3$ are defined by
\begin{align*}
u_1\,=\,z_1,\quad u_2\,=\,z_2z_3,\quad u_3\,=\,z_2+z_3.
\end{align*}
The relation $g(z_1,z_2,z_3)\,=\,0\,$ in (\ref{relzch3}) becomes
\begin{align}
u_1+u_1^2-u_3\,=\,0.\nonumber
\end{align}
We conclude that $k(x_1,x_2)^{\langle\sigma\rangle}\,=\,k(u_1,u_2,u_3)\,=\,k(u_1,u_2)$ is
rational over $k$.
The generators $u_1,u_2$ of $k(x_1,x_2)^{\langle\sigma\rangle}$ over $k$ is given
in terms of $x_1,x_2$ as follows:
\[
u_1\,=\, \frac{x_1}{x_2^2},\quad u_2\,=\,\frac{-(x_1^2+x_1x_2^2+x_2^3)}{x_2^5}.
\]

\bigskip


\end{document}